\documentclass[a4paper,12pt]{article} 

\usepackage{amsmath,amssymb, mleftright}      
\numberwithin{equation}{section}

\newtheorem{thm}{Theorem}[section]

\newtheorem{lem}{Lemma}[section]

\newcommand{\pt}{\partial}

\newcommand{\re}{\mathbb R}
\newcommand{\Nt}{\mathbb N}

\newcommand{\al}{\alpha}

\newcommand{\ep}{\varepsilon}
\newcommand{\lam}{\lambda}

\renewcommand{\t}{\tau}

\def\<{\langle }

\title{Improvement of the general theory \\ for one dimensional nonlinear wave equations \\ related to the combined effect}
\author{
Shu Takamatsu
\footnote{
Master course, Mathematical Institute,
Tohoku University,
Aoba, Sendai 980-8578, Japan.
e\_mail: shu.takamatsu.r8@dc.tohoku.ac.jp }
}
\date{{\small\it }
\[
\begin{array}{ll}
\mbox{\footnotesize{\bf Keywords:}}
& \mbox{\footnotesize nonlinear wave equation, one dimension,
lifespan.}\\
\mbox{\footnotesize{\bf MSC2020:}}
& \mbox{\footnotesize primary 35L70, secondary 35A01}\\
\end{array}
\]
}
\begin{document}
\maketitle

%%%%%%%%%%%%%%%%%%%%%%%%%%%%%%%%%%%%%%%%%%%%%%%
%%%%%%%%%%%%%%%%%%% ABSTRACT %%%%%%%%%%%%%%%%%%%%%
%%%%%%%%%%%%%%%%%%%%%%%%%%%%%%%%%%%%%%%%%%%%%%% 

\begin{abstract}
We focus on the general theory to the Cauchy problem for one dimensional nonlinear wave equations with small initial data.
In the general theory, we aim to obtain the lower bound estimate of  the lifespan of classical solution.
In this paper, we improve it in some case related to the \lq\lq combined effect'', which was expected complete more than 30 years ago.
\end{abstract}

%%%%%%%%%%%%%%%%%%%%%%%%%%%%%%%%%%%%%%%%%%%%%%%%%%
%%%%%%%%%%%%%%%%%%%%% SECTION1 %%%%%%%%%%%%%%%%%%%%%%%
%%%%%%%%%%%%%%%%%%%%%%%%%%%%%%%%%%%%%%%%%%%%%%%%%%

\section{Introduction}
Let us consider the Cauchy problem for nonlinear wave equations;
\begin{equation}
\label{IVP}
\left\{
\begin{array}{ll}
	 u_{tt}-u_{xx}=F(u, Du, Du_x)
	&\mbox{in}\quad \re\times(0,T),\\
	u(x,0)=\ep f(x),\ u_t(x,0)=\ep g(x)
	& x\in\re,
\end{array}
\right.
\end{equation}
where 
$
D=(\pt_t, \pt_x),
$
$f, g\in C^{\infty}_0(\re)$, $T>0$, and $\ep>0$ is a small palameter.
Let
$
\lam \in \re^5.
$
We assume that in a neighborhood of $\lam=0$, $F(\lam)$ is a $C^{\infty}$ function satisfying
\begin{equation}
\label{order_F}
F(\lam)=\mathcal{O}(|\lam|^{1+\al}),
\end{equation}
where $\al\geq1$ is an integer.
This paper is aimed at studying the lifespan $T(\ep)$, the maximal existence time of classical solutions of (\ref{IVP}).

By virtue of Li, Yu and Zhou \cite{LYZ91, LYZ92}, the result of the general theory to one dimensional nonlinear wave equations is as follows.
\begin{equation}
\label{general_theory}
T(\ep)\geq
\left\{
\begin{array}{lll}
C\ep^{-\al/2} & \mbox{in general},\\
C\ep^{-\al(\al+1)/(\al+2)} & \mbox{if}\  \displaystyle\int_{\re}g(x)dx=0, \\
C\ep^{-\min\{\beta_0/2,\al\}} & \mbox{if}\ \pt^{\beta}_uF(0,0,0)=0,\ \al+1 \leq \forall\beta \leq \beta_0,
\end{array}
\right.
\end{equation}
where $\beta_0\geq\al+1$ is an integer.
However, Morisawa, Sasaki and Takamura \cite{MST}, and Kido, Sasaki, Takamatsu and Takamura \cite{KSTT}
implied that there is a possibility to improve the general theory by showing the so-called \lq\lq combined effect".
For instance, according to a part of \cite{KSTT}, we can obtain the following lower bound estimate
for nonlinear term $u^pu_t^q+u^r$ with integer $p,q,r>1$;
\begin{equation}
\label{combined}
T(\ep)\geq C\ep^{-(p+q)(r-1)/(r+1)}
\end{equation}
in the case of
\begin{equation}
\label{power}
\frac{r+1}{2}<p+q<r
\end{equation}
and
\begin{equation}
\label{condition_g}
\int_{\re}g(x)dx=0.
\end{equation}
On the other hand, (\ref{general_theory}) produces in this case that
\begin{equation}
\label{general}
T(\ep)\geq C\ep^{-\max\{(p+q-1)(p+q)/(p+q+1),\ (r-1)/2\}}.
\end{equation}
Because of the condition (\ref{power}), we can verify that the estimate (\ref{combined}) is better than the estimate (\ref{general}). 
Consequently, by setting $\al = p+q-1$ and $\beta_0=r-1$, we can conjecture improvement of the general theory
if (\ref{condition_g}) and
\begin{equation}
\label{condition_F}
\pt^{\beta}_uF(0,0,0)=0 \ (\al+1 \leq \forall\beta \leq \beta_0)\quad \mbox{with}\quad  \al+1 \leq \beta_0 <2\al
\end{equation}
is fulfilled.

Then, our result is the following;
\begin{equation}
\label{result}
T(\ep) \geq C \ep^{-(\al+1)\beta_0/(\beta_0+2)}
\end{equation}
in the case of (\ref{condition_g}) and (\ref{condition_F}). This result implies (\ref{combined}).
Moreover, when (\ref{condition_g}) and (\ref{condition_F}) hold, the general theory (\ref{general_theory}) yields
\begin{equation}
\label{a}
T(\ep)\geq C\ep^{-\max\{\al(\al+1)/(\al+2),\ \beta_0/2\}}.
\end{equation}
Because 
\[
\displaystyle \frac{(\al+1)\beta_0}{\beta_0+2}\geq\displaystyle \frac{\beta_0}{2},\  \frac{\al(\al+1)}{\al+2}
\]
follows from  $\al+1\leq\beta_0<2\al$,
we can verify our result (\ref{result}) is better than (\ref{a}).

Summing up, the general theory is improved as follows;
\[
T(\ep)\geq
 \left\{
 \begin{split}
 &
  \begin{array}{lll}
   C\ep^{-\al/2} & \mbox{in general},\\
   C\ep^{-\al(\al+1)/(\al+2)} & \mbox{if}\ \displaystyle\int_{\re}g(x)dx=0, \\
  \end{array} \\
 &
  \begin{array}{lll} 
   C\ep^{-\min\{\beta_0/2, \al\}}  &\mbox{if}\ \pt^{\beta}_uF(0,0,0)=0\ \mbox{for}\ \al+1 \leq \forall\beta \leq \beta_0,\\
   C\ep^{-(\al+1)\beta_0/(\beta_0+2)}& \mbox{if}\ \pt^{\beta}_uF(0,0,0)=0\ \mbox{for}\ \al+1 \leq \forall\beta \leq \beta_0<2\al\\
  \end{array}  \\
 &\qquad\qquad\qquad \mbox{and} \displaystyle\int_{\re}g(x)dx=0.\\
 \end{split}
 \right. \\
\]

Finaly, we mention the sharpness of our result (\ref{result}).
When there exists $\beta_0$ such that $\al+1\leq \beta_0 < 2\al$, 
we consider nonlinear term $|u_t|^{\al+1}+|u|^{\beta_0+1}$ or $|u|^p|u_t|^q+|u|^{\beta_0+1}$ with $p,q\in \Nt_{>1}$ and $p+q=\al+1$.
Then, the blow up part of \cite{MST} or \cite{KSTT} yield the sharpness of our result (\ref{result}).

%%%%%%%%%%%%%%%%%%%%%%%%%%%%%%%%%%%%%%%%%%%%%%%%%%
%%%%%%%%%%%%%%%%%%%%% SECTION2 %%%%%%%%%%%%%%%%%%%%%%%
%%%%%%%%%%%%%%%%%%%%%%%%%%%%%%%%%%%%%%%%%%%%%%%%%%

\section{Preliminaries and Main result}

We note that to investigate Cauchy problem (\ref{IVP}), it suffices, essentialy, to consider the following Cauchy problem of quasi-linear hyperbolic equations;
\begin{equation}
\label{IVP_quasi}
\left\{
\begin{array}{ll}
	 u_{tt}-u_{xx}=b(u,Du)u_{xx}+2a_0(u,Du)u_{tx}+ F(u,Du)
	&\mbox{in}\quad \re\times(0,T),\\
	u(x,0)=\ep f(x),\ u_t(x,0)=\ep g(x)
	& x\in\re.
\end{array}
\right.
\end{equation}
For example, see Section 7 in T.Li, and Y.Zhou \cite{Zhou_book} for this fact.
Let $\hat{\lam}\in \re^3$. From the assumption (\ref{order_F}), we assume that 
\begin{equation}
\label{order_b}
b(\hat{\lam}), a_0(\hat{\lam})=\mathcal{O}(|\hat{\lam}|^{\al}),
\end{equation}
\begin{equation}
F(\hat{\lam})=\mathcal{O}(|\hat{\lam}|^{1+\al}).
\end{equation}
In addition, the assumption (\ref{condition_F}) corresponds to
\begin{equation}
\label{condi_F}
\pt^{\beta}_uF(0,0)=0 \ (\al+1 \leq \forall\beta \leq \beta_0),
\end{equation}
where  $\al+1 \leq \beta_0 <2\al$. Let us set
\[
F(v,Dv)=F(v,0)+\tilde{F}(v,Dv)Dv.
\]
Then by assumption (\ref{condi_F}), we have $F(\tilde{\lam},0)=\mathcal{O}(|\tilde{\lam}|^{1+\beta_0})$ and $\tilde{F}(\hat{\lam})=\mathcal{O}(|\hat{\lam}|^{\al})$,
where $\tilde{\lam}\in\re$. Hereafter, we consider the Cauchy problem (\ref{IVP_quasi}).

Throughout this paper, we assume that the initial data $f$ and $g$ satisfies
\begin{equation}
\label{supp_initial}
\mbox{\rm supp }f,\ \mbox{supp }g\subset\{x\in\re:|x|\le R\},\quad R\geq1.
\end{equation}
Let $u$ be a classical solution of (\ref{IVP_quasi}) in the time interval $[0,T]$.
Then the support condition of the initial data, (\ref{supp_initial}), implies that
\begin{equation}
\label{support_sol}
\mbox{supp}\ u(x,t)\subset\{(x,t)\in\re\times[0,T]:|x|\leq t+R\}.
\end{equation}
For example, see Appendix of John \cite{John_book} for this fact.
Now, we divide $\mbox{supp}\ u(x,t)$ into two domains, $D$ and $D'$, where $D$ is the interior domain;
\[
D:=\{(x,t)\in\re\times[0,T]\ :\ t-|x|\geq R\},
\]
and $D'$ is the exterior domaim;
\[
D'=D^{c}\cap\{(x,t)\in\re\times[0,T]\ :\ |x|\leq t+R\}.
\]
It is well-known that $u$ satisfies the following integral equation.
\begin{equation}
\label{u}
u(x,t)=\ep u^0(x,t)+L(b(u,Du)u_{xx}+2a_0(u,Du)u_{tx}+ F(u,Du))(x,t),
\end{equation}
where  a linear integral operator $L$ for a function $v=v(x,t)$ in Duhamel's term is defined by
\begin{equation}
\label{nonlinear}
L(v)(x,t):=\frac{1}{2}\int_0^tds\int_{x-t+s}^{x+t-s}v(y,s)dy,
\end{equation}
and $u^0$ is a solution of the free wave equation with the same initial data,
\begin{equation}
\label{u^0}
u^0(x,t):=\frac{1}{2}\{f(x+t)+f(x-t)\}+\frac{1}{2}\int_{x-t}^{x+t}g(y)dy.
\end{equation}
When g satisfies the condition (\ref{condition_g}), the strong Huygens' principle
\begin{equation}
\label{Huygens}
u^0(x,t)\equiv0\quad\mbox{in}\ D
\end{equation}
holds.

By the Sobolev embedding theorem, there exists $E_0>0$ so small that 
\[
\|f\|_{L^{\infty}}\leq1 ,\ \forall f\in H^1(\re), \|f\|_{H^1(\re)}<E_0
\]

For any given integer $S\geq4$ and any given suitably small positive number $E$ and T, we introduce the following set of function;
\[
X_{S,E,T}:=\{v(x,t): D_{S,T}(v)\leq E, \pt_t^lv(x,0)=u_l^{(0)}(x)\ (l=0,1,\ldots, S)\},
\]
where $u_0^{(0)}(x):=\ep f(x)$, $u_1^{(0)}(x):=\ep g(x)$ and $u_l^{(0)}(x)$ is the value of  $\pt_t^lu(x,t)$ at $t=0$, and $D_{S,T}(v)$ is following;
\[
D_{S,T}(v):=\|v\|_1+\|v\|_2+ \sup_{ 0\leq t \leq T}\|Dv(\cdot,t)\|_{D,S,2},
\]
where
\[
\|v\|_1=\sup_{(x,t)\in D}(1+t+x)^{p}|v(x,t)| \quad  \mbox{with}\quad  p=\frac{2\al-\beta_0}{\beta_0(\al+1)}>0 ,  
\]
\[
\|v\|_2=\sup_{(x,t)\in  D'}|v(x,t)| ,
\]
\[
\|Dv(\cdot,t)\|_{D,S,2}=\sum_{|k|\leq S}\|D^kDv(\cdot,t)\|_{L^2({\re})}.
\]
Moreover, let $\tilde{X}_{S,E,T}$ be the subset of $X_{S,E,T}$ composed of all elements in $X_{S,E,T}$ with compact support (\ref{support_sol}) in the variable $x$ for any fixed $t\in [0,T]$.
 
The main result is following.
\begin{thm}
\label{thm}
Assume (\ref{condition_g}) and (\ref{order_b})-(\ref{condi_F}).
For any given integer $S\geq4$, there exist positive constants $\ep_0$ and $C_0$ with $C_0\ep\leq E_0$ such that the Cauchy problem (\ref{IVP_quasi}) admits
a unique classical solution $u\in\tilde{X}_{S,C_0\ep,T(\ep)}$ on $[0,T(\ep)]$ for any $0<\ep<\ep_0$.
Then, $T(\ep)$ satisfies
\[
T(\ep) \geq c \ep^{-(\al+1)\beta_0/(\beta_0+2)}
\]
where $c$ is a positive constant.

Moreover, with eventual modification on a set with zero measure in the variable $t$, we have
\[
u\in C([0,T(\ep)]; H^{S+1}(\re)),
\]
\[
u_t\in C([0,T(\ep)]; H^{S}(\re)),
\]
\[
u_{tt}\in C([0,T(\ep)]; H^{S-1}(\re)).
\]
\end{thm}

%%%%%%%%%%%%%%%%%%%%%%%%%%%%%%%%%%%%%%%%%%%%%%%%%%
%%%%%%%%%%%%%%%%%%%%% SECTION3 %%%%%%%%%%%%%%%%%%%%%%%
%%%%%%%%%%%%%%%%%%%%%%%%%%%%%%%%%%%%%%%%%%%%%%%%%%

\section{The proof of Theorem \ref{thm}}

In order to prove theorem \ref{thm}, we define a map
\[
M:v \rightarrow u=Mv
\]
by solving the following Cauchy Problem for linear wave equations;
\[
\left\{
\begin{array}{ll}
	 u_{tt}-u_{xx}=b(v,Dv)u_{xx}+2a_0(v,Dv)u_{tx}+ F(v,Dv)
	&\mbox{in}\quad \re\times(0,T),\\
	u(x,0)=\ep f(x),\ u_t(x,0)=\ep g(x)
	& x\in\re.
\end{array}
\right.
\]

We note that the outline of the proof of Theorem \ref{thm} is same as the one of Li, Yu and Zhou \cite{LYZ91, LYZ92}.
Lemma 3.4 and Lemma 3.5 in Li, Yu and Zhou \cite{LYZ91, LYZ92} is only replaced with following lemmas, respectively.

\begin{lem}
\label{dst_u}
When $E>0$ is suitably small, for any given $v \in \tilde{X}_{S,E,T}$, $u=Mv$ satisfies
\[
D_{S,T}(u)\leq C_1\{\ep+(R(E,T)+\sqrt{R(E,T)})(E+D_{S,T}(u))\},
\] 
where $C_1$ is a positive constant independent of $E$ and $T$, and
\begin{equation}
\label{RET}
R(E,T)=E^{\al}(1+T)^{1+p}+E^{\beta_0}(1+T)^{2-\beta_0p}.
\end{equation}
\end{lem}

\par\noindent
{\bf Proof.} See the section4 below.
\hfill$\Box$

\begin{lem}
\label{dst_u-}
Under the assumptions of Lemma \ref{dst_u}, for any given $v,\ v' \in \tilde{X}_{S,E,T}$, if $u=Mv,\ u'=Mv'$ satisfy $u, u'\in \tilde{X}_{S,E,T}$, then we have
\[
D_{S-1,T}(u-u')\leq C_2(R(E,T)+\sqrt{R(E,T)})(D_{S-1,T}(u-u')+D_{S-1,T}(v-v')),
\] 
where $C_2$ is a positive constant independent of $E$ and $T$, and $R(E,T)$ is still defined by (\ref{RET}).
\end{lem}

\par\noindent
{\bf Proof.} See Lemma 3.5 in Li, Yu and Zhou \cite{LYZ92} and prove this in similar way to the proof of Lemma \ref{dst_u}.
\hfill$\Box$

\vspace{3mm}
Now, we start to prove Theorem \ref{thm}.
Take $C_0=3\max\{C_1, C_2\}$, where $C_1,\ C_2$ are positive constants appearing in Lemma \ref{dst_u} and \ref{dst_u-}, respectively.
Moreover, we take $E=E_\ep=C_0\ep$ and $T=T_\ep=c\ep^{-(\al+1)\beta_0/(\beta_0+2)}-1$. Then, noting that (\ref{RET}) and the definition of $p$,
\[
R(E_\ep,T_\ep)+\sqrt{R(E_\ep,T_\ep)}=C_0^{\al}c^{1+p}+C_0^{\beta_0}c^{2-\beta_0p}+\{C_0^{\al}c^{1+p}+C_0^{\beta_0}c^{2-\beta_0p}\}^{1/2}.
\]
Therefore, we can choose a suitable constant $c$, such that 
\[
R(E_\ep,T_\ep)+\sqrt{R(E_\ep,T_\ep)}\leq1/C_0.
\]
Using Lemma \ref{dst_u} and \ref{dst_u-}, for any given $v, v' \in \tilde{X}_{S,E_\ep,T_\ep}$, $u=Mv, u'=Mv'$ satisfy
\[
D_{S,T_\ep}(u)\leq E_\ep
\]
and
\[
D_{S-1,T_\ep}(u-u')\leq \frac{1}{2}D_{S-1,T_\ep}(v-v').
\]
Consequently, similarly to  Li, Yu and Zhou \cite{LYZ91, LYZ92}, we can verify the map $M$ possesses a unique fixed point. 
The proof of Theorem \ref{thm} is now finished.

%%%%%%%%%%%%%%%%%%%%%%%%%%%%%%%%%%%%%%%%%%%%%%%%%%
%%%%%%%%%%%%%%%%%%%%% SECTION4 %%%%%%%%%%%%%%%%%%%%%%%
%%%%%%%%%%%%%%%%%%%%%%%%%%%%%%%%%%%%%%%%%%%%%%%%%%

\section{The proof of Lemma \ref{dst_u}}
First, we note next lemma.
\begin{lem}
For any $v\in X_{S,E,T}$, we have
\begin{equation}
\label{H1}
\sup_{ 0\leq t \leq T}\sum_{0\leq |k|\leq S}\|D^kv(\cdot,t)\|_{L^{\infty}({\re})}\leq CE,
\end{equation}
where $C$ is positive constant.
\end{lem}
\par\noindent
{\bf Proof.} Noting the definition of $X_{S,E,T}$, from the Sobolev embedding theorem we know
\[
H^1(\re) \subset L^{\infty}(\re)
\]
is a continuous embedding, so (\ref{H1}) holds.
\hfill$\Box$

\vspace{3mm}
Here a positive constant $C$ independent of $E$, $T$ and $\ep$ may change from line to line.

\vskip10pt
\par\noindent
{\bf The estimate of $\bf \|u\|_1$.}
\par

First, we estimate $\|u\|_1$. 
We assume $(x,t) \in D$. (\ref{u}) and (\ref{Huygens}) yield that
\begin{equation}
\label{u+}
\begin{split}
u(x,t)=&L(b(v,Dv)u_{xx}+2a_0(v,Dv)u_{tx}+ \tilde{F}(v,Dv)Dv)(x,t)\\ 
& +L(F(v,0))(x,t).
\end{split}
\end{equation}
Noting that the definition of $\|\cdot\|_1$, the second term has following estimate.
\[
\begin{split}
|L(F(v,0))(x,t)|
&\leq |L(\chi_DF(v,0))(x,t)|+|L(\chi_{D'}F(v,0))(x,t)| \\
&\leq CE^{\beta_0+1}\int_0^tds\int_{x-t+s}^{x+t-s}(1+y+s)^{-(\beta_0+1)p}\chi_D(y,s)dy\\
&\quad +CE^{\beta_0+1}\int_0^tds\int_{x-t+s}^{x+t-s}\chi_{D'}(y,s)dy \\
&\leq CE^{\beta_0+1}\int_0^tds\int_{x-t+s}^{x+t-s}(1+y+s)^{-(\beta_0+1)p}\chi_D(y,s)dy\\
&\quad +CE^{\beta_0+1}(1+T), \\
\end{split}
\] 
where $\chi_{D}$ is a characteristic function of $D$. We employ the change of variables $\xi=s+y,\ \eta=s-y$.
Then, noting that 
\begin{equation}
\label{positiveness}
-(\beta_0+1)p+1=\frac{(\beta_0-\al)(2+\beta_0)}{(\al+1)\beta_0}>0,
\end{equation}
we have
\[
\begin{split}
&\int_0^tds\int_{x-t+s}^{x+t-s}(1+y+s)^{-(\beta_0+1)p}\chi_D(y,s)dy\\
&=\int^{t-x}_Rd\eta\int^{x+t}_R(1+\xi)^{-(\beta_0+1)p}d\xi \\
& \leq C(t-x-R)(1+t+x)^{-(\beta_0+1)p+1} \\
& \leq C(1+T)^{-(\beta_0+1)p+2}.
\end{split}
\]
The inequality $1+T \leq (1+T)^{-(\beta_0+1)p+2}$ holds by (\ref{positiveness}), so that we obtain
\[
|L(F(v,0))(x,t)| \leq CE^{\beta_0+1}(1+T)^{-(\beta_0+1)p+2}.
\]

Next, we estimate the first term in (\ref{u+}). Integration by parts yields that
\[
\begin{split}
&|L(b(v,Dv)u_{xx})(x,t)|=C\left|\int_0^tds\int_{x-t+s}^{x+t-s}b(v,Dv)u_{xx}(y,s)dy\right| \\
&\leq C \left|\int_0^t b(v,Dv)u_{xx}(x+t-s,s) - b(v,Dv)u_{xx}(x-t+s,s) ds\right| \\
&\quad + C \left|\int_0^tds\int_{x-t+s}^{x+t-s}\frac{\pt b(v,Dv)(y,s)}{\pt y}u_{x}(y,s)dy\right| \\
&\leq CE^{\al}D_{S,T}(u)(1+T)+CE^{\al} \int_0^tds\int_{x-t+s}^{x+t-s}|Dv(y,s)||u_{x}(y,s)|dy \\
&\leq CE^{\al}D_{S,T}(u)(1+T)+CE^{\al}\int_0^t\|Dv(\cdot,s)\|_{L^2(\re)}\|u_{x}(\cdot,s)\|_{L^2(\re)}ds \\
&\leq CE^{\al}D_{S,T}(u)(1+T).
\end{split}
\]
Similarly, we have
\[
|L(a_0(v,Dv)u_{xx})(x,t)|\leq CE^{\al}D_{S,T}(u)(1+T)
\]
and
\[
|L(\tilde{F}(v,Dv)\pt_x v)(x,t)| \leq CE^{\al+1}(1+T).
\]
Regarding $L(\tilde{F}(v,Dv)\pt_t v)(x,t)$, by inverting the order on $(y,s)$-integral, we have similar estimate;
\[
|L(\tilde{F}(v,Dv)\pt_t v)(x,t)| \leq CE^{\al+1}(1+T).
\]
Consequentry, we have
\begin{equation}
\label{esti_1}
\|u\|_1=\sup_{(x,t)\in D}(1+t+x)^{p}|u(x,t)| \leq CR(E,T)(E+D_{S,T}(u)).
\end{equation}

%%%%%%%%%%%%%%%%%%%%%%%%%%%%%%%%%%%%%%%%%%%%%%%%%%%%%%%%%%%%%%%%%%%%%%%%%%%%%%%%%%%%%%%%%%%%%%%%%%%%
\vskip10pt
\par\noindent
{\bf The estimate of $\bf \|u\|_2$.}
\par
We assume $(x,t) \in D'$. Moreover we consider $x\geq 0$.
Then, we have
\[
\begin{split}
&|L(b(v,Dv)u_{xx})(x,t)|\leq CE^{\al}D_{S,T}(u) \int_0^tds\int_{x-t+s}^{x+t-s}dy\\
&\leq CE^{\al}D_{S,T}(u) \int_0^tds\int_{s-R}^{s+R}dy\\
&\leq CE^{\al}D_{S,T}(u)(1+T)
\end{split}
\]
Similarly, we have
\[
|L(a_0(v,Dv)u_{xx})(x,t)|\leq CE^{\al}D_{S,T}(u)(1+T)
\]
and
\[
|L(F(v,Dv))(x,t)|\leq CE^{\al+1}(1+T).
\]
The case of $x \leq 0$ is similar to the one above, so we omit the details. Therefore, we have
\[\begin{split}
|u(x,t)|
&\leq\ep |u^0(x,t)|+|L(b(u,Du)u_{xx}+2a_0(u,Du)u_{tx}+ F(u,Du))(x,t)| \\
&\leq C\{\ep+ E^{\al}(1+T)(D_{S,T}(u)+E) \} \\
&\leq C\{\ep+ R(E,T)(D_{S,T}(u)+E) \}
\end{split}\]
and
\begin{equation}
\label{esti_2}
\|u\|_2\leq C\{\ep+ R(E,T)(D_{S,T}(u)+E) \}.
\end{equation}

%%%%%%%%%%%%%%%%%%%%%%%%%%%%%%%%%%%%%%%%%%%%%%%%%%%%%%%%%%%%%%%%%%%%%%%%%%%%%%%%%%%%%%%%%%%%%%%%%%%%
\vskip10pt
\par\noindent
{\bf The estimate of $\bf \|Du(\cdot,t)\|_{D,S,2}$.}
\par
For any double index $k=(k_1,k_2)$ with $0\leq |k|\leq S$, we obtain the following energy integral formula;
\[
\begin{split}
&\|D^ku_t(\cdot, t)\|^2_{L^2(\re)}+ \int_{\re}a(v,Dv)(x,t)(D^ku_x(x,t))^2 dx\\
&= \|D^ku_t(\cdot, 0)\|^2_{L^2(\re)}+ \int_{\re}a(v,Dv)(x,0)(D^ku_x(x,0))^2 dx\\
&+\int^t_0\int_{\re}\frac{\pt b(v,Dv)(x, \t)}{\pt \t}(D^ku_x(x,\t))^2 dxd\t \\
&-2\int^t_0\int_{\re}\frac{\pt b(v,Dv)(x, \t)}{\pt x}(D^ku_x(x,\t))(D^ku_\t(x,\t)) dxd\t \\
&-2\int^t_0\int_{\re}\frac{\pt a_0(v,Dv)(x, \t)}{\pt x}(D^ku_\t(x,\t))^2 dxd\t \\
&+2\int^t_0\int_{\re}G_k(x, \t)D^ku_\t(x,\t) dxd\t  +2\int^t_0\int_{\re}g_k(x, \t)D^ku_\t(x,\t) dxd\t \\
&=: \|D^ku_t(\cdot, 0)\|^2_{L^2(\re)}+ \int_{\re}a(v,Dv)(x,0)(D^ku_x(x,0))^2 dx \\
&+I_1+I_2+I_3+I_4+I_5,
\end{split}
\]
where $a(\cdot)=1+b(\cdot)$,
\[
G_k=D^k(b(v,Dv)u_{xx})-b(v,Dv)D^ku_{xx}+2(D^k(a_0(v,Dv)u_{xx})-a_0(v,Dv)D^ku_{xx}),
\]
\[
g_k=D^kF(v,Dv).
\]
For example, see Section 8 in T.Li, and Y.Zhou \cite{Zhou_book} for this fact. (\ref{order_b}) and (\ref{H1}) yield that
\[\begin{split}
|I_1|
&\leq CE^{\al}\int^t_0\int_{\re}(D^ku_x(x,\t))^2 dxd\t \\
&\leq CE^{\al}\int^t_0\|D^ku_x(\cdot,\t)\|_{L^2(\re)} d\t \\
&\leq CE^{\al}D^2_{S,T}(u)(1+T).
\end{split}\]
Similarly, we have
\[
|I_2|, |I_3|\leq CE^{\al}D^2_{S,T}(u)(1+T).
\]

Noting that $0\leq |k|\leq S$, we have
\[
\|D^k(b(v,Dv)u_{xx})-b(v,Dv)D^ku_{xx}\|_{L^2(\re)} \leq CE^{\al}D_{S,T}(u).
\]
For the term of involving $a_0$, we have similar estimate. Hence,
\[
|I_4|\leq CE^{\al}D^2_{S,T}(u)(1+T)
\]
holds.

Now, we estimate $|I_5|$.
When $|k|\geq 1$, we have
\[
\begin{split}
|I_5|
&\leq CE^{\al}\int^t_0\int_{\re}|DvD^ku_\t(x,\t)| dxd\t \\
&\leq CE^{\al}\int^t_0\|Dv(\cdot,\t) \|_{L^2(\re)}\|D^ku_\t(\cdot,\t) \|_{L^2(\re)}d\t \\
&\leq CE^{\al+1}D_{S,T}(u)(1+T). \\
\end{split}
\]
When $k=0$, we have
\[
\begin{split}
I_5
&=2\int^t_0\int_{\re}F(v,Dv)(x, \t)u_\t(x,\t) dxd\t \\
&=2\int^t_0\int_{\re}F(v,0)(x, \t)u_\t(x,\t) dxd\t +2\int^t_0\int_{\re}\tilde{F}(v,Dv)Dv(x, \t)u_\t(x,\t) dxd\t \\
&=:J_1+J_2
\end{split}
\]
The estimate of $J_2$ is similar to the one above, so we have
\[
|J_2|\leq CE^{\al+1}D_{S,T}(u)(1+T).
\]  
In terms of $J_1$, we can calculate as following.
\[
\begin{split}
J_1
&=2\int^t_0\int^{\t+R}_{-\t-R}F(v,0)(x, \t)u_\t(x,\t) dxd\t \\
&=2\left(\int^{-R}_{-t-R}\int^{t}_{-x-R}+\int^{R}_{-R}\int^{t}_{0}+\int^{t+R}_{R}\int^{t}_{x-R} \right)F(v,0)(x, \t)u_\t(x,\t) d\t dx  \\
&=:K_1+K_2+K_3.
\end{split}
\]
\[
\begin{split}
|K_1| 
&\leq 2\left|\int^{-R}_{-t-R}\{F(v,0)(x, t)u(x,t)-F(v,0)(x, -x-R)u(x,-x-R)\}dx\right| \\
&\quad +2\left|\int^{-R}_{-t-R}\int^{t}_{-x-R} \frac{\pt F(v,0)(x, \t)}{\pt \t}u_\t(x,\t) d\t dx\right|\\
&\leq 2E^{1+\beta_0}D_{S,T}(u)(1+T) + 2E^{\beta_0}\left|\int^{-R}_{-t-R}\int^{t}_{-x-R} v_\t(x,\t) u_\t(x,\t) d\t dx\right|\\
&\leq CE^{1+\beta_0}D_{S,T}(u)(1+T). 
\end{split}
\]
Similarly, we have
\[
|K_2|, |K_3| \leq CE^{1+\beta_0}D_{S,T}(u)(1+T).
\]
Therefore, in the case of $k=0$, 
\[
|I_5|\leq C(E^{\al}+E^{\beta_0})ED_{S,T}(u)(1+T)
\]
holds.

Consequently, we have
\begin{equation}
\label{esti_3}
\|Du(\cdot,t)\|_{D,S,2}\leq C\{\ep+\sqrt{R(E,T)}(E+D_{S,T}(u))\}.
\end{equation}

%%%%%%%%%%%%%%%%%%%%%%%%%%%%%%%%%%%%%%%%%%%%%%%%%%%%%%%%%%%%%%%%%%%%%%%%%%%%%%%%%%%%%%%%%%%%%%%%%%%%
\vskip10pt
\par\noindent
{\bf The proof of Lemma \ref{dst_u}.}
\par
Summing up (\ref{esti_1}), (\ref{esti_2}), and (\ref{esti_3}), the proof of Lemma \ref{dst_u} is now completed.

%%%%%%%%%%%%%i%%%%%%%%%%%%%%%%%%%%%%%
%%%%%%%%%%%% Acknowledgement %%%%%%%%%%%%%%%
%%%%%%%%%%%%%%%%%%%%%%%%%%%%%%%%%%%%%
\section*{Acknowledgement}
\par
The author would like to thank Professor Hiroyuki Takamura (Tohoku Univ.,
Japan) for introducing me this problem, his consistent encouragements and discusstions.

%%%%%%%%%%%%%%%%%%%%%%%%%%%%%%%%%%%%%%
%%%%%%%%%%%% References %%%%%%%%%%%%%%%%%%%%
%%%%%%%%%%%%%%%%%%%%%%%%%%%%%%%%%%%%%%

\bibliographystyle{plain}

\end{document}